\documentclass[]{amsart}
\usepackage{amsmath,paralist,amssymb}
\usepackage{soul,color,graphics}
\newif\ifdeveloping
\developingtrue

\usepackage[ocgcolorlinks]{hyperref}

\ifdeveloping

\fi

\DeclareMathOperator{\Hom}{Hom}
\DeclareMathOperator{\Ker}{Ker}
\DeclareMathOperator{\Sym}{Sym}
\DeclareMathOperator{\Alt}{Alt}
\DeclareMathOperator{\Epi}{Epi}
\DeclareMathOperator{\Aut}{Aut}
\DeclareMathOperator{\Out}{Out}
\DeclareMathOperator{\Inn}{Inn}
\newcommand{\PRA}{\operatorname{PRA}}
\newcommand{\PSL}{\operatorname{PSL}}
\newcommand{\PSO}{\operatorname{PSO}}
\newcommand{\dist}{\operatorname{dist}}

\newcommand{\C}{\mathbb{C}}
\newcommand{\F}{\mathbb{F}}
\newcommand{\bbH}{\mathbb{H}}
\newcommand{\N}{\mathbb{N}}
\newcommand{\R}{\mathbb{R}}
\newcommand{\Z}{\mathbb{Z}}
\newcommand{\Q}{\mathbb{Q}}
\newcommand{\calB}{\mathcal{B}}
\newcommand{\calD}{\mathcal{D}}
\newcommand{\calH}{\mathcal{H}}
\newcommand{\calS}{\mathcal{S}}
\newcommand{\calP}{\mathcal{P}}
\newcommand{\calE}{\mathcal{E}}
\newcommand{\la}{\langle}
\newcommand{\twola}{\la\!\la}
\newcommand{\ra}{\rangle}
\newcommand{\twora}{\ra\!\ra}
\newcommand{\ulin}{\underline}
\newcommand{\olin}{\bar}
\newcommand{\oline}{\overline}
\newcommand{\varpbt}{\varphi_B((\underline{t}))}
\newcommand{\varpt}[1]{\varphi_{#1}((\underline{t}))}

\numberwithin{equation}{section}
\makeatletter
\let\c@equation\c@subsection
\makeatother

\newcommand{\thra}{\twoheadrightarrow}

\newcommand{\calL}{\mathcal{L}}

\newcommand{\be}[1]{\begin{equation}\label{#1}}
\newcommand{\ee}{\end{equation}}

\theoremstyle{plain}

\newtheorem*{thm*}{Theorem}
\newtheorem{subthm}{Theorem}[subsection]

\newtheorem*{prop*}{Proposition}

\newtheorem*{claim*}{Claim}
\newtheorem*{conj*}{Conjecture}
\newtheorem{conj}{Conjecture}[section]
\newtheorem{thms}[conj]{Theorem}
\newtheorem{lems}[conj]{Lemma}
\newtheorem{props}[conj]{Proposition}
\newtheorem{claim}[conj]{Claim}
\newtheorem{prob}[conj]{Problem}

\newtheorem{cors}[conj]{Corollary}

\theoremstyle{definition}
\newtheorem{defns}[conj]{Definition}

\newtheorem{ques}[conj]{Question}
\newtheorem*{que*}{Question}
\newtheorem{sque}{Question}

\newcounter{subtheorem}
\newenvironment{subtheorem}{\refstepcounter{conj}}{}

\author{Alexander Lubotzky}
\address{Einstein Institute of Mathematics,
Hebrew University,
Jerusalem 91904
Israel
}
\email{alexlu@math.huji.ac.il}

\title[Dynamics of $\Aut(F_n)$ actions]{Dynamics of $\Aut(F_n)$ actions\\ on group presentations and representations}

\begin{document}

\ifdeveloping
\fontsize{11}{16}\selectfont
\fi

\newcommand{\ml}[1]{\colorbox{yellow}{$#1$}}

\dedicatory{
\begin{flushright}
\normalsize{To Bob Zimmer}
\end{flushright}
}
\begin{abstract}
Several different areas of group theory, topology and geometry have led to the study of the action of $\Aut(F_n)$ --- the automorphism group of the free group on $n$ generators --- on $\Hom(F_n,G)$ when $G$ is either finite, compact or simple Lie group.
In this survey, we describe these topics and results, with special emphasis on some similarities and with an effort to give a somewhat uniform treatment.
This perspective sometimes suggests new questions, conjectures and methods  borrowed from one area to another.
\end{abstract}

\maketitle

\tableofcontents

\section{Introduction}\label{section1}
In this paper we will survey a collection of results originated from various different contexts but all have the following common form:

Let $F=F_n$ be the free group on $n\ge 2$ generators and $G$ a group.
Denote
\begin{align}
H_n(G)&=\Hom(F_n,G)\label{eq:1.1}\\
E_n(G)&=\Epi(F_n,G)=\{\varphi\in H_n(G)|\varphi \text{ is onto}\} \label{eq:1.2}
\end{align}
The group $\Aut(G)\times\Aut(F_n)$ acts on $H_n(G)$ by:
\begin{multline}\label{eq:1.3}
(\beta,\alpha)(\varphi)=\beta\circ\varphi\circ\alpha^{-1}\\
 \text{where }\alpha\in\Aut(F_n),\,\beta\in\Aut(G)\text{ and }\varphi\in H_n(G).
\end{multline}
The action clearly preserves $E_n(G)$.
We denote by $\bar{H_n}(G)$ (resp.\ $\bar{E_n}(G)$) the quotient set $\Aut(G)\backslash H_n(G)$ (resp.\ $\Aut(G)\backslash E_n(G)$).
So $\Aut (F_n)$ acts on  $\bar H_n(G)$ and on $\bar E_n(G)$.
What really acts on $\bar E_n(G)$ is $\Out(F_n)$, as $\Inn(F_n)$ acts trivially.
The orbits of $\bar E_n(G)$  under $\Aut (F_n)$ are called ``$T$-systems of $G$'' in some of the old literature  about the subject, but we will not use this term.
What we will present here is a collection of results and methods to study these sets and these actions.
The motivations come from various quite different areas of research: presentation of groups, actions of groups on handlebodies, computational group theory and the product replacement algorithm, the theory of linear groups, Fuchsian and Kleinian groups, compact groups and more.

Here is a brief outline of the content of this survey, chapter by chapter:

The set $\bar E_n(G)/\Aut(F_n) $ is in a one to one correspondence with equivalence classes of presentations of $G$ using $n$ generators.
While for an infinite groups $G$, this can be a very large and complicated set, it is known in some cases, and conjectured in others, that it is a ``tame set'' when $G$ is finite.
In \S\ref{section2}, we describe some of the known results and open problems.
The orbits of $\Aut(F_n)$ on $E_n(G)$ can be illustrated as the connected components of a graph --- the PRA graph.

The product replacement algorithm is a probabilistic algorithm providing a pseudo random element in a finite group given by a set of generators.
It is really a random walk on the PRA graph whose vertices are the $n$-generating sets, i.e., the set $E_n(G)$.

Furthermore, The vertices of the graph and its connected
components are also in correspondence with actions of $G$ on
handlebodies and their equivalence classes. So, the topological
question and the computational group theory problem discussed in
\S\ref{section3}, are different forms of the questions raised in
\S\ref{section2} on the presentation theory of $G$.

In \S\ref{section4}, we take the opportunity to present in passing another graph associated with $G$ - the Andrews--Curtis graph.
Its origin is the classical Andrews--Curtis conjecture in combinatorial group theory and topology, but it also got a new interest from computational group theory.
In Sections \ref{section5} and \ref{section6} respectively, we discuss the connectivity of the $\PRA$ graphs for solvable and simple group, respectively.
The case of finite simple groups is of special interest.
A long standing conjecture (Wiegold conjecture) suggests various extensions to compact and Lie groups as well as potential applications to the representation theory of $\Aut(F_n)$.

In Section \ref{section7}, we will compare the situation to the case where $F$ is replaced by its profinite completion $\hat F$ and explain why profinite presentations ``behave nicer'' than discrete presentations:
All is due to a beautiful lemma of Gash\"utz, which seems to be not as well-known as it should.

In Section \ref{section8}, we will treat the case where $G$ is a semisimple compact Lie group.
Here $\Hom(F_n,G)$ and $\Epi(F_n,G)$ are the same form  a measure theoretical point of view (where by $\Epi(F_n,G)$ we mean now, the homomorphisms $\varphi$ for which $\varphi(F_n)$ is dense in $G$).
We study the ergodicity of the action of $\Aut(F_n)$ on $H_n(G)=\Hom(F_n,G)=G^n$.

Let now $G$ be a general semisimple  Lie group.
Here almost nothing is known, except for the cases $G=\PSL_2(\R)$ and $\PSL(\C)$.
We present in Section \ref{section9} some basic questions and few results.

In the last chapter we will replace $F_n$ by $T_g=\pi_1(S_g)$ the
fundamental group of a closed surface of genus $g$. The case of
surface groups deserves a separate survey but we will touch it
only briefly, suggesting along the way a possible program for
proving that the mapping class groups are not linear. The group
$\Aut (T_g)$ acts on $\Hom(T_g,G)$ and $\Out(T_g)$ on the
equivalence classes mod $\Inn(G)$. Note that $\Out(T_g)$ is the
mapping class group $M_g$ and by restricting ourselves to a subset
of $\Hom(T_g,\PSL_2(\R))$ we recover the action of $M_g$ on the
Fricke--Teichm\"uller  space. The case of $G$ compact has also
been studied in the literature and even the case of $G$ finite
came up recently in the work of Dunfield and Thurston on finite
covers of random 3-manifolds.

These notes are based on a series of lectures given at Yale
University in February 2008. I would like to thank the
participants for their remarks and questions and, in particular,
to Yair Minsky for information about Kleinian Groups. I am also
grateful to Nir Avni, Shelly Garion, Tsachik Gelander and Yair
Glasner for several discussions.

This paper is dedicated to Bob Zimmer from whom I learned to search for the ``big picture''.
Bob has always been a friend as well as a source of inspiration in his leadership as a scientist and as an administrator.
%
%

\section{Presentations}\label{section2}
Let $G$ be a finitely generated  group.
By $d(G)$ we denote the minimal number of generators of $G$.
A presentation of $G$ is an exact sequence
\be{eq:2.1a}
1\to R\to F_n\xrightarrow{\varphi}G\to1.
\ee
where $F=F_n$ is the free group on $n$ generators $x_1,\ldots,x_n$.
Clearly $n\geq d(G)$.
We denote by $d_F(R)$ the minimal number of generators of $R$ as a {\it normal} subgroup of $F$.

A basic question in the theory of presentations of groups is to what extend a presentation as in \eqref{eq:2.1a} is unique.
Clearly we can ``twist'' \eqref{eq:2.1a} by the action of $\Aut(G)\times \Aut(F_n)$ in the following way
\be{eq:2.1b}
(\beta,\alpha)(\varphi)=\beta\circ \varphi\circ\alpha^{-1}\quad \text{where $(\beta,\alpha)\in \Aut(G)\times \Aut(F_n)$}.
\ee
Another way to change a presentation is by enlarging $n$.
Let $\pi:F_{n+1}\to F_n$ be the natural epimorphism sending $x_i$ (as an element of $F_{n+1}$) to $x_i$ (as an element of $F_n$) for $i=1,\ldots,n$ and $x_{n+1}$ to the identity.
Then
\be{eq:2.1c}
1\to N\to F_{n+1} \xrightarrow{\varphi\circ\pi}G\to 1
\ee
is also a presentation of $G$, which we call a lifting of \eqref{eq:2.1a}.

\subsection{Waldhausen's problem:}\label{subsec2.2}
Waldhausen (see \cite[p.~92]{LySc}) raised the question whether every presentation of a group $G$ with $n>d(G)$ can be obtained from a presentation with $n=d(G)$ by a sequence of liftings.
An equivalent formulation is:
\begin{que*}
Let $F_m$ be the free group on $m$ generators and $N\lhd F_m$.
Assume $d(G)< m$ where $G=F_m/N$.
Does $N$ contain a primitive element of $F_m$?
\end{que*}

Recall that $y\in F_m$ is primitive if it is part of some basis of $F_m$ or equivalently, it is in the orbit of $x_1$ under the action of $\Aut(F_m)$.
\subsection{Gruenberg's Questions}
It was probably B.\ H.\ Neuman who already in the 1930's considered in a systematic way the connection between various presentations of the same group $G$.
Let us follow Gruenberg's treatment in \cite{Gru} who presented three fundamental questions:
\begin{sque}\label{question2.3a}
Given two presentations
\[
1\to R_i\to F_n \xrightarrow{\varphi_i}G\to 1
\]
$i=1,2$.
Is $d_{F_n}(R_1)=d_{F_n}(R_2)$?
\end{sque}
\begin{sque}\label{question2.3b}
Let $r(G)$ be the minimum number of relations needed to define $G$, i.e., the minimum of $d_{F_n}(R)$ over all possible presentations of $G$ as in \eqref{eq:2.1a}.
Is $r(G)$ realized in a minimal presentation (i.e.\ one in which $n=d(F_n)=d(G))$?
\end{sque}
\begin{sque}\label{question2.3c}
Is $d_{F_n}(R)-d(F_n)$ independent of the presentation \eqref{eq:2.1a} and therefore an invariant of $G$?
\end{sque}
Of course, a positive answer to Question \ref{question2.3c} would imply a similar one for \ref{question2.3b} and \ref{question2.3a}.

Unfortunately, the answer to these questions in general is
negative. For example, Dunwoody and Pietrowski \cite{DP} showed
that the group $G=\langle a,b: a^2=b^3\rangle$ is a one relator
group which has also a presentation with 2 generators which needs
more than one relater (see also \cite{Gru} and references
therein).

Another example, attributed to G.\ Higiman is given in \cite[p.\ 93]{LySc} where it is shown that the Baumslag--Solitar group $G=\langle x,y : x^{-1}y^2 x=y^3\rangle$ is also generated by $x$ and $z=y^4$, but has no presentation with these two generators with a single defining relation.

This shows that the answer to Question \ref{question2.3a} and \ref{question2.3c} is negative.
Noskov \cite{Nos} (see also Evans \cite{Eva93A}) showed that the answer to Waldenhaus' question \ref{subsec2.2} is negative.
We do not know a counter example to Question \ref{question2.3b}, but one very likely exists.

One can easily see that if for every $n$, the action of $\Aut(G)\times\Aut(F_n)$ on the set of possible presentations of $G$ (with $n$ generators) is transitive, then the answer to the first three questions raised in this chapter would be positive.
But this is far from being the case if $G$ is infinite.
The case of finite $G$ is more delicate  and more interesting, and this will be the subject of the next sections.
Before elaborating on it, we will introduce the convenient language of the product replacement graphs.

\section{The Product Replacement Algorithm and its graph}\label{section3}

\subsection{}\label{subsection3.1}
  The {\it Product Replacement Algorithm} {\it (PRA)} is a practical algorithm for generating random elements of a finite group.
The algorithm was introduced and analyzed in \cite{CLMNO}.
Although it has no rigorous justification, practical experiments have shown excellent performance.
It quickly became a popular algorithm for generating random group elements, and was included in two frequently used group computation packages: GAP and MAGMA.

The Product Replacement Algorithm can be described as a random walk on a graph, called the {\it Product Replacement Graph} (or the {\it PRA Graph}).
It will be more convenient for us to look at the following extended graph.
For any $n \geq d(G)$, let
\[
V_n(G) = \{(g_1,\ldots, g_n)\in G^n : \langle g_1,\ldots,
g_n\rangle = G\}
\]
be the set of all {\it generating $n$-tuples} of $G$.

The {\it extended PRA graph}, denoted $\tilde{X}_n(G)$, has $V_n(G)$ as its set of vertices.
The edges correspond to the following so-called {\it Nielsen moves} $R^{\pm}_{i,j}, L^{\pm}_{i,j} , P_{i,j} , I_i$ for $1\leq i \neq j \leq n$, where
\begin{align*}
R^{\pm}_{i,j} &: (g_1,\ldots,g_i,\ldots,g_n) \to (g_1,\ldots,g_i\cdot g^{\pm1}_j,\ldots,g_n)\\
L^{\pm}_{i,j} &: (g_1,\ldots,g_i,\ldots,g_n) \to (g_1,\ldots,g_j^{\pm1}\cdot g_i,\ldots,g_n)\\
P_{i,j} &: (g_1,\ldots,g_i,\ldots,g_j,\ldots,g_n) \to (g_1,\ldots,g_j,\ldots,g_i,\ldots,g_n)\\
I_i & : (g_1,\ldots,g_i,\ldots,g_n) \to (g_1,\ldots,g^{-1}_i,\ldots,g_n)
\end{align*}

Strictly speaking, the Product Replacement Algorithm is a random walk on a subgraph $X_n(G)$ of $\tilde{X}_n(G)$, which is obtained by removing the edges corresponding to the $P_{i,j}$ and $I_i$.
The output of the algorithm is a random entry chosen from the tuple at the end of the random walk.
As observed in \cite[Proposition 2.2.1]{Pa}, when $n \geq d(G) + 1$, the graph $X_n(G)$ is connected if and only if $\tilde{X}_n(G)$ is connected.
The connected components of $\tilde X_n(G)$ are also called Nielsen equivalence classes of $n$-generating sets of $G$.

\subsection{}\label{subsection3.2}
Recall now the well-known result of Nielsen (cf.\ \cite{LySc})
that $\Aut(F_n)$ is generated by the Nielsen moves
$\{R^{\pm}_{i,j},L^{\pm}_{i,j},P_{i,j},I_i\}_{1\leq i,j\leq n}$
viewed as automorphisms of $F_n$ acting on the $n$-tuple of
generators $(x_1,\ldots,x_n)$. Moreover, $V_n(G)$ is naturally
identified with $E_n(G)$ the set of epimorphisms from $F_n$, the
free group on $x_1,\ldots,x_n$, onto $G$, since every such
epimorphism $\varphi$ is uniquely defined by the generating vector
$(\varphi(x_1),\ldots,\varphi(x_n))$. It follows, therefore, that
the connected components of the graph $\tilde{X}_n(G)$ are exactly
the orbits of $\Aut(F_n)$ acting on $E_n(G)$. This observation
plays a crucial role in \cite{LP} when the mixing rate of the
random walk on $\tilde{X}_n(G)$ is related to the possibility of
$\Aut (F_n)$ (or some variants of it) have property $(T)$ or
$(\tau)$. We will not elaborate on this issue here --- referring
the reader to \cite{LP} for more information and for some
interesting open problems.

We denote by $\bar X_n(G)$ the quotient graph $\Aut(G)\backslash \tilde X_n(G)$.
Its vertices are in one to one correspondence with the points of $\bar E_n (G)$.
These in turn are in one to one correspondence with the set of normal subgroups $N$ of $F_n$ with $F_n/N$ isomorphic to $G$.
The action of $\Aut(F_n)$ on $\olin E_n(G)=\olin X_n(G)$ factor through $\Out(F_n)=\Aut(F_n)/\Inn(F_n)$.

The connectivity of $\tilde{X}_n(G)$ (resp., $\bar X_n(G)$) is equivalent to $|E_n(G)/\Aut(F_n)|=1$ (resp., $|\bar E_n(G)/\Aut(F_n)|=1$) and if $\tilde X_n(G)$ is connected so is $\bar X_n(G)$.
But we do not know if the converse is true in general.
In \cite[Proposition 2.4.1]{Pa} it is observed that this is the case if $n\geq2d(G)$.
Anyway, the recent interest in the product replacement algorithm put forward the question of the connectivity of $\tilde{X}_n(G)$ which is essentially equivalent to the question we started with, in Section 2, i.e., whether $G$ has an essentially unique presentation on $n$ generators.
>From now on we will use both languages.

\smallskip
There are very few general results in this context which holds for every finite group.
Here is one:

\begin{props}\label{proposition3.1}
Let $G$ be a finite group and $\mu(G)$ the maximal size of a minimal set of generators of $G$ (i.e.~if $m>\mu(G)$, every set $S$ of generators of $G$ with $|S|\geq m$ contains a proper subset of generators).
Then $\tilde X_n(G)$ is connected for every $n\geq \mu(G)+d(G)$.
\end{props}

\begin{proof}
Let $g_1,\ldots,g_d$ be a set of generators of $G$, $d=d(G)$ and $(\ulin g)=(1,\ldots,1,g_1,\ldots,g_d)\in\tilde X_n(G)$.
Let $(\ulin x)=(x_1,\ldots,x_n)\in \tilde X_n(G)$ be an arbitrary vector.
By assumption it contains a subset $S$ of size $\mu=\mu(G)$ of generators.
After changing order we can assume that these are $x_1\ldots,x_{\mu}$.
We can now move $(\ulin x)$ a series of Nielsen moves to $(x_1,\ldots, x_{\mu},1,\ldots,1,g_1,\ldots,g_d)$ and the latter can be moved to $(\ulin g)$. Thus $\tilde X_n(G)$ is connected.
\end{proof}

Very few results seem to be known on $\mu(G)$: Whiston \cite{Wh}
and Saxl--Winston \cite{WS} estimate it for $S_n$ and $PSL_2(q)$,
respectively. Nikolov (unpublished) showed that for every finite
simple group of Lie type $G$ of rank at most $r$ over  a finite
field of order at most $p^e$, $\mu(G)\leq f(r,e)$ where $f$
depends on $r$ and $e$ but not on $p$.

Anyway, for every finite group $G$, one can easily see that $d(G)$ and $\mu(G)$ are both bounded by $\log _2(|G|)$.
So we can deduce that $\tilde X_n(G)$ is connected for any $n\geq 2 \log_2(|G|)$.
But as observed in \cite{MW} one can do even better:

Let $l=l(G)$ be the maximum length of a chain of strictly decreasing non-trivial subgroups of $G$.
It is easy to see that $d(G)\leq\mu(G)\leq l(G)$.

\begin{props}\label{newprop3.2}
If $n>l(G)$ then $\tilde X_n(G)$ in connected.
\end{props}

\begin{proof}
Fix a generating vector $(s_1,\ldots,s_d)$ with $d=d(G)$.
Let $(\ulin t)=(t_1,\ldots,t_n)\in \tilde X_n (G)$.
We will show that $(\ulin{t})$ is connected to $(s_1,\ldots,s_d,1,\ldots,1)$.

Put $G_i=\la t_1,\ldots, t_i\ra$, so $G=G_n\geq G_{n-1} \geq \cdots \geq G_1\geq \{1\}$.
Since $n>l(G)$, $G_j=G_{j-1}$ for some $j>0$, so $t_j$ is a word in $t_1,\ldots,t_{j-1}$.
So $(\ulin t)$ is connected to  $(t_1,\ldots,t_{j-1},1,t_{j+1},\ldots, t_n)$.
We can (by changing names) assume that $t_1$ of $(\ulin t)$ satisfies $t_1=1$.
Thus $G=\la t_2,\ldots, t_n\ra$ and hence $(\ulin t)$ is connected to $(s_1,t_2,\ldots,t_n)$.
Define now $G_1=\la s_1\ra$ and $G_i=\la s_1,t_2,\ldots,t_i\ra$ for $i\geq 2$.
Again, as $n>l(G)$, we must have $G_j=G_{j-1}$ for some $j$ and since $s_1\neq 1$, $j>1$.
So $t_j\in\la s_1,t_2,\ldots,t_{j-1}\ra$ and $(\ulin t)$ is equivalent to $(s_1,t_2,\ldots,t_{j-1},1,t_{j+1},\ldots,t_n)$ and hence to $(s_1,s_2,t_3,\ldots,t_m)$.
We continue by induction to deduce that $(\ulin t)$ is connected to $(s_1,\ldots,s_d,t_{d+1},\ldots,t_n)$.
As $\la s_1,\ldots,s_d\ra = G$, it implies that $(\ulin t)$ is connected to $(s_1,\ldots s_d,1,\ldots,1)$ as promised.
\end{proof}

Since $l(m)\leq \log _2(|G|)$ we deduce:

\begin{cors}\label{newcorollary3.3}
$\tilde X_n(G)$ is connected for $n>\log_2(|G|)$.
\end{cors}

\subsection{Free action of finite groups on handlebodies}\label{subsection3.3}
The Neilsen equivalence classes of generating sets of $G$ parametrize free actions of the group $G$ on handlebodies as we will now explain, following \cite{MW} and the references therein.

Let $\calH$ be an orientable  three dimensional handlebody of genus $g\geq 1$.
Two (effective) actions $\rho_1,\rho_2: G \to \text{Homeo}(\calH)$ are said to be \emph{equivalent} if there is a homeomorphism $h:\calH\to \calH$ such that $h\rho_1(g)h^{-1}=\rho_2(g)$  for each $g\in G$.
They are \emph{weakly equivalent} if their images are conjugate, i.e., there is $h\in \text{Homeo}(\calH)$ such that $h\rho_1(G)h^{-1}=\rho_2(G)$.
Equivalently, there is $\alpha\in \Aut (G)$ such that $h\rho_1(g)h^{-1}=\rho_2(\alpha(g))$, i.e., $\rho_1$ and $\rho_2\circ\alpha$ are equivalent.

>From now on when we talk about actions of $G$ on $\calH$ we mean orientation preserving $\emph{free}$ actions.
We will assume $g\geq 1$, as the only free action on the handlebody of genus $0$, the 3-ball, is by the trivial group.

If $G$ acts on $\calH$ freely, then the quotient map $\calH \to \calH/G$ is a covering map which induces an extension
\be{equation3.3.1}
1\to \pi_1(\calH)\to \pi_1(\calH/G)\to G\to 1.
\ee
Note that $\pi_1(\calH)$ is a free group on $g$ generators and by \cite[Theorem 5.2]{He} $\calH/G$ is also a handlebody, so $\pi_1(\calH/G)$ is also a free group on, say, $n$ generators.
Nielsen--Schreier Theorem implies that $g=1+|G|(n-1)$ i.e.\ $n=1+\frac1{|G|}(g-1)$.
Conversely, if we start with a handlebody $\calH'$ of genus $n$, so $\pi_1(\calH')=F_n$, every epimorphism onto $G$ gives rise to a covering $\calH$ which is a handlebody of genus $g$, on which $G$ acts.
This sets up a surjective map from the family of (free) actions of $G$ on the handlebody $\calH$ of genus $g$ onto the set of epimorphisms from $F_n$ onto $G$.
The latter is exactly the set of vertices of $\tilde X_n(G)$.
Now, every automorphism of $F_n=\pi_1(\calH/G)$ is induced by an homeomorphism of $\calH/G$.
The following theorem it is deduced in \cite[Theorem 2.3]{MW} using elementary arguments from covering theory:

\begin{thms}
Let $G$ be a finite group, $n\geq 1$ and $g=1+|G|(n-1)$.
The equivalence classes of actions of $G$ on a genus $g$ handlebody correspond bijectively to the Nielsen equivalence classes of $n$-generating sets of $G$, i.e., to the connected components of $\tilde X_n(G)$.
The weak equivalence classes of these actions correspond to the connected components of $\olin X_n(G)$.
\end{thms}

The theorem shows that all the results discussed in this survey on the connected components of $\tilde X_n(G)$ or $\olin X_n(G)$ have direct topological applications and in fact are equivalent to such topological statements.
The correspondence goes even further.
An important notion in the study of actions of $G$ on $\calH$ is ``stabilization'': two actions $\rho_1$ and $\rho_1$ can be equivalent after ``adding'' one (or more) handle to $\calH$.
This is equivalent to the question whether the two corresponding $n$-generating set $(\ulin t)=(t_1,\ldots,t_n)$ and $(\ulin s)=(s_1,\ldots,s_n)$ are on the same connected component when considered as $(s_1,\ldots,s_n,1)$ and $(t_1,\ldots,t_n,1)$ in $\tilde X_{n+1}(G)$.
So, all these seemingly pure algebraic questions, to be discussed later, carry a significant amount of topological information.
We will usually stick to the algebraic language, leaving the reader the translation to this topological setting.

Later in the paper we will show various connectivity results.
It is a highly non-trivial problem to show that two $n$-generating sets are \emph{not} Nielsen equivalent (except for $n=2$ see Section \ref{subsection5.2} below).
Some powerful methods using Fox calculus were developed in \cite{Lustig} and \cite{LM93} but as far as we know these methods have not been applied as of now to yield non-equivalence for $n$-generating sets of finite groups.

\section{The Andrews--Curtis conjecture and its graph}\label{section4}
Before going to a more detailed study of $\tilde{X}_n(G)$, let us mention in passing another graph associated with a group --- the Andrews--Curtis graph $AC_n(G)$.
This graph also has its roots in a deep problem in topology and in presentation theory, but the interest in it has revived recently from the point of view of computational group theory.

Let $G$ be a group generated by a finite symmetric set $h_1,\ldots,h_d$, and $N$ a normal subgroup of $G$.
As usual $d_G(N)$ denotes the minimal number of elements of $N$ generating $N$ as a normal subgroup of $G$.
For $n\geq d_G(N)$ we define the graph $AC_n(G,N)$ as follows:

Its vertices are the $n$-tuples $(g_1,\ldots,g_n)\in N^n$ with $\langle\langle g_1,\ldots, g_n\rangle\rangle_G=N$, i.e., those which generate $N$ as a normal subgroup of $G$.

A vector $(g_1,\ldots,g_n)$ will be connected to its image under the moves $L_{i,j}^{\pm}$, $R_{i,j}^{\pm}$, $P_{i,j}$, $I_i$ as in \ref{subsection3.1}
as well as $(g_1,\ldots,g_n)\to(g_1,\ldots,h_jg_ih_j^{-1},\ldots,g_n)$ for every $1\leq i\leq n$ and $1\leq j\leq d$.

Note that all moves indeed take an $n$-tuple of normal generating set of $N$ to another one.
The case when $N=G$ is of special interest.
In this case we write $AC_n(G)$ instead of $AC_n(G,G)$.

The famous Andrew--Curtis conjecture is equivalent to:
\begin{conj}[Andrews--Curtis \cite{AC}]
The graph $AC_n(F_n)$ is connected.
\end{conj}

The conjecture is usually expressed in a different language: Note that a vector of $AC_n(F_n)$ amounts to a vector of $n$ elements of the free group $F_n$, i.e., $n$ words in $x_1,\ldots,x_n$ which normally generate $F_n$.
In other words, this is  a presentation of the trivial group by $n$ generators and $n$ relations.
The Andrews--Curtis conjecture predicts that any such presentation is obtained from the standard presentation $\la x_1,\ldots,x_n;x_1,\ldots,x_n\ra$ using a finite series of Nielsen moves or conjugations.
This is one of the most outstanding conjectures in combinatorial group theory, with various potential applications to topology (see \cite{AC}).

Let us mention that in general $AC_n(G)$ is not connected.
E.g.\ if $G=\F_p^n$ (or any abelian group with $d(G)=n$) then $AC_n(G)$ is the same as $\tilde{X}_n (G)$ and, as will be shown in the next section (in a different language!), the number of the connected component of $\tilde{X}_n(G)$ is $(p-1)/2$ and the graph is not connected.

The topic has received a new interest in recent years from computational group theory \cite{BKM}:
A natural generalization of the $\PRA$ algorithm is the following algorithm to produce a pseudo random element in the normal closure $N$ of a given set of elements $g_1,\ldots,g_n$ inside a group $G$ generated by given generators $h_1,\ldots,h_d$.
The algorithm starts with the vector $(g_1,\ldots,g_n)$ of $AC_n(G,N)$ and takes a random walk on the graph.
The output is a random component from the end vector of the random walk.

Many researchers believe that the Andrews--Curtis conjecture is false and some (see  \cite{BKM} and the references therein) have tried to disprove it by using computer calculations.
The idea is that if $\pi:H\to G$ is an epimorphism of groups then $\pi$ induces a graph theoretic map
$\tilde{\pi}: AC_n(H)\to AC_n(G)$ (not necessarily onto!).
If one could find one finite group $G$ and an epimorphism $\pi:F_n\to G$ such that $\tilde{\pi}(AC_n(F_n))$ is not connected then the conjecture is false.
Various calculations have been performed till it was shown in \cite{BLM}.
\begin{thms}\label{theorem3.3}
Let $G$ be a finite group and $n\geq\max\{d_G(G),2\}$.
Then two vectors of $AC_n(G)$ are connected by a path iff their images in $AC_n(G/[G,G])$ are connected.
\end{thms}

The connected components of the Andrews--Curtis graph of abelian group are easy to understand and one can deduce:

\begin{cors}\label{corollary3.4}
For every epimorphism $\pi:F_n\to G$ where $G$ is a finite group, the image $\tilde{\pi}(AC_n(F_n))$ is connected.
\end{cors}

This ``finitary Andrews--Curtis Conjecture'' does not give much insight on the original conjecture, but shows that the computational efforts carried out in order to disprove it, have all been in vain.

We end up this section by giving a sketch of the proof for a special case of the Theorem --- the case when $G$ is perfect.
In this case $G/[G,G]=\{1\}$ and the Theorem claims that $AC_n(G)$ is connected.
Let's prove it:

First, denote by $M(G)$ the intersection of all maximal normal subgroups of $G$.
An easy observation is that a subset $\{y_1,\ldots,y_k\}$ of $G$ normally generate $G$ if and only if it generates it $\bmod  M(G)$.
For infinite groups, this is usually a useless observation (e.g., $M(F)=\{e\}$ if $F$ is a free group).
But for a finite group $G$, $G/M(G)$ is always a direct product of finite simple groups, and if $G$  is also perfect, all the simple groups are non-abelian.
We can replace $G$ by $G/M(G)=\prod_{i=1}^r S_i$, $S_i$ non-abelian finite simple groups and every $z=(z_1,\ldots,z_r)\in\prod_{i=1}^r S_i$ with $z_i\neq e$ for all $1\leq i\leq r$, generates $G/M(G)$ normally.
To  show now that $AC_n(G)$ is connected, we show that every vector in $AC_n(G)\subseteq G^n$ is connected to $(z,1,\ldots,1)$, where $z$ is an  element as above.

So let $(y_1,\ldots,y_n)\in AC_n(G)$. Look at $y_{n-1}$ and $y_n$
as elements of $G/M(G)=\prod_{i=1}^{r} S_i$. We can conjugate
$y_{n-1}$ by some $g\in G$ s.t.\ $\tilde y_{n-1}={y}_{n-1}^g y_n$
is not the identity in every component, unless \emph{both}
$y_{n-1}$ and $y_n$ are identities at that component: to do so
simply conjugate such that the conjugation of $y_{n-1}$ in
component $i$ is different from the component of $y_n$  there.
Then  use the normal closure of the new $\tilde{y}_{n-1}$, which
is the product of all components in which $\tilde y_{n-1}$ is
non-identity, to ``clean'' $y_n$, i.e.\ to make it the identity in
these components and hence altogether to replace $y_n$ by the
identity.

One can continue like that to get that $(y_1,\ldots,y_n)$ is connected  to $(z',1,\ldots,1)\in G^n$ where all components of $z'$ as an element of $\prod_{i=1}^r S_i$ are non-identity.
As $n\geq 2$, we can change $(z',1,\ldots,1)$ to $(z',z,1,\ldots,1)$, than switch to $(z,z',1,\ldots,1)$ and then ``clean'' $z'$ by $z$ to get $(z,1,\ldots,1)$ as promised.

\section{Finite solvable groups}\label{section5}
 While the theory of presentations of infinite groups seems wild and one expects very few general results to hold, the theory of presentations of finite groups may have some pleasant properties.
 In particular, the question we are interested in, i.e., the number of orbits of $\Aut(G)\times\Aut(F_n)$ acting on the $n$-generators presentations of $G$, is quite well-understood for finite solvable groups due to the following two results of Dunwoody.

\begin{subtheorem}
\begin{subthm}[\cite{Du2}]\label{thm4.1a}
Let $G$ be a finite solvable group with $d(G)<n$.
Then $\Aut(F_n)$ acts transitively on $E_n(G)$, i.e., $\tilde{X}_n(G)$ (and hence also $\olin X_n(G)$) is connected.
\end{subthm}

The assumption $d(G)<n$ is crucial.
Note that for $G=\F_p^n$ the $n$-dimensional vector space over the field $\F_p$, $E_n(G)$ is the set of all bases of $G$, so can be identified with $GL_n(\F_p)$.
The action of $\Aut(F_n)$ on it is via $SL_n^{\pm}(\F_p)=\{A\in GL_n(\F_p)|\det(A)=\pm 1\}$ and so $|E_n(G)/\Aut(F_n)|=(p-1)/2$.
On the other hand $\Aut(G)=GL_n(\F_p)$ and so $|\bar{E}_n(G)/\Aut(F_n)|=1$.
But there are groups where even the second set is large:

\begin{subthm}[\cite{Du1}]
For every $2\leq n\in \N$, every $k\in\N$ and every prime $p$, there exists a finite $p$-group of nilpotency class two with $|\bar{E}_n(G)/\Aut(F_n)|\geq k$.
In particular, the $\PRA$ graph  $\tilde{X}_n(G)$ has at least $k$ components.
\end{subthm}
\end{subtheorem}

Theorem \ref{thm4.1a} shows that for finite solvable groups, Questions \ref{subsec2.2} and \ref{question2.3b} have affirmative answer.
The same is true for \ref{question2.3a} provided $n>d(G)$.

The proof of Theorem \ref{thm4.1a} is non-trivial but quite elementary.
As this is the main positive result known, let us sketch its proof:

Let $\{e\}=M_0\lhd M_1\lhd\cdots\lhd M_r=G$ be a chief series of $G$, i.e., $M_i\lhd G$ and $M_{i+1}/M_i$ is a minimal normal subgroup of $G/M_i$.
It is, therefore, an irreducible $\F_p[G/M_{i+1}]$ module, for some prime $p$.
Let $h_1,\ldots,h_{n-1}$ be a finite set of generators for $G$ and $(g_1,\ldots,g_n)\in V_n(G)$.
We will argue by induction on $r$, so assume the theorem for $r-1$.
This means that we can move by Nielsen transformations from $(g_1,\ldots,g_n)$ to $(m,m_1h_1,\ldots,m_{n-1}h_{n-1})$ where $m,m_1,\ldots,m_{n-1}\in M_1$.
We should show that the last vector is connected to $(e,h_1,\ldots,h_{n-1})$.

We can assume $m\neq e$.
Indeed, if $m=e$, then $m_1h_1,\ldots,m_{n-1}h_{n-1}$ generate $G$ (since $(e,m_1h_1,\ldots,m_{n-1}h_{n-1})\in V_n(G)$) and so we can create any word of them in the first component.

So assume that $m\neq e$.
Think of $M_1$ as an additive group, an $G/M_1$-module.
Note also that as $M_1$ is abelian, the action of $m_ih_i$ on $M_1$ by conjugation is the same as that of $h_i$.
As $M_1$ is irreducible, one deduces that every $m'\in M_1$ can be written as a sum of conjugates of $m$ by words in $m_1h_1,\ldots, m_{n-1}h_{n-1}$ and this is so for every $m\neq e$.
Now for each $i=1,\ldots,n-1$, write $m_i$ as a sum $m_i=m_{i1}+\cdots+m_{it}$ where $m_{ij}$ is a conjugate of $m$.
Then we can change $m$ to $m_{ij}$ to clean $m_{ij}$ out of $m_i$ and gradually clean $m_i$.
Do it for all $i=1,\ldots,n-1$ to get $(\tilde{m},h_1,\ldots,h_{n-1})$ and then eliminate $\tilde{m}$, which is possible since $h_1,\ldots,h_{n-1}$ generate $G$.

\section{Finite simple groups}\label{section6}
\subsection{}
We are coming to the most interesting case, with some very interesting open problems and some potential applications --- the case when $G$ is a finite non-abelian simple group.
For such a group $G$ the leading long-standing conjecture is the following one, which is attributed to Jim Wiegold:

\begin{conj}\label{conjecture5.1}
Let $G$ be a finite simple group.
Then $\olin X_n(G)$ is connected for every $n\geq 3$.
\end{conj}

Note, that by the classification of the finite simple groups, it is known that for such $G$, $d(G)\leq 2$, so the conjecture combined with Dunwoody's theorem \ref{thm4.1a} made Pak \cite{Pa} to ask:

\begin{ques}\label{question5.2}
Is it true that for every finite group $G$, $\olin X_n(G)$ is connected if $n>d(G)$?
\end{ques}

It is interesting to mention that Dunwoody in a review on \cite{Gi} (see Math.\ Review MR0435226 in 1997) wrote:
``It seems unlikely that this result is true for an arbitrary finite group $G$''.
But in the years since then no counter example has been found, so maybe the answer to Question \ref{question5.2} is indeed positive.
As an intermediate step one can suggest the following conjecture which looks more feasible, as it is known to be true for simple groups (see Theorem \ref{thms5.5} below) and for solvable groups (Theorem \ref{thm4.1a} above).

\begin{conj}\label{conj5.3}
Let $G$ be a finite group, $(\ulin t)=(t_1,\ldots,t_n)$ and $(\ulin s)=(s_1,\ldots,s_n)\in \tilde X_n(G)$.
Then $(t_1,\ldots,t_n,1)$ and $(s_1,\ldots,s_n,1)$ are connected in $\olin X_{n+1}(G)$.
\end{conj}

Note that it is easy to see that $(\ulin s)$ and $(\ulin t)$ become connected to each other in $\tilde X_{n+d}(G)$ when $d=d(G)$ (see \cite[Proposition 6.1]{MW}).
Note also that Conjecture \ref{conj5.3} has a topological equivalent formulation, as hinted in Section \ref{subsection3.3}, it asserts that any two (free) actions of $G$ on a handlebody become equivalent after adding one handle (see \cite{MW} for more on stabilizations of actions).

We will come again to the case where $n>d(G)$, but let us first clear up the situation when $n=d(G)$.

\subsection{$G$ simple and $n=2$}\label{subsection5.2}
When $G$ is a finite simple group and $n=d(G)=2$, the situation is very much different than what is predicted by Conjecture \ref{conjecture5.1}.

\begin{thms}[Garion-Shalev \cite{GS08}]\label{theorem5.3}
Let $G$ be a non-abelian finite simple group. Then $|\bar{E}_2(G)/\Aut(F_2)|\to \infty$ when $|G|\to \infty$, or equivalently the number of connected components of $\bar X_2(G)$ is going to infinity with $G$.
\end{thms}

The special case of $G=\PSL_2(q)$ was proved by Guralnick and Pak \cite{GP} who also conjectured the general case.

Let us sketch the proof: A classical result of Nielsen asserts
that if $\alpha\in\Aut F_2$, with $F_2=F(x,y)$ the free group in
$x$ and $y$, then the commutator $[\alpha(x),\alpha(y)]$ is
conjugate in $F_2$ to either $[x,y]$ or $[x,y]^{-1}$. This implies
that if $\varphi,\psi\in\Epi(F_2,G)$ are on the same
$\Aut(F_2)\times\Aut(G)$ orbits then $[\varphi(x),\varphi(y)]$ is
conjugate to $[\psi(x),\psi(y)]^{\pm 1}$ in $\Aut(G)$. So to prove
the theorem, Garion and Shalev showed that ``almost all the
elements of $G$ are commutators of pairs of generators of $G$'',
i.e.\ the proportion of these elements in $G$ goes to one when the
order of $G$ goes to infinity. Once this is proved, it follows
that the number of components of $\olin X_2(G)$ grows at least as
fast as the number of the conjugacy classes of $G$. The claim
above is proved by combining two methods. First, note that the
function $F(g)=\#\{(a,b)\in G\times G| [a,b]=g\}$ is a class
function on $G$, i.e., constant on conjugacy classes. So, by
harmonic analysis on finite groups, it can be expressed using
characters of G. Moreover, a classical result of Frobenius from
1896 gives an explicit formula:
\[
F(g)=|G|\cdot\sum_{\chi}\frac{\chi(g)}{\chi(1)}.
\]
Thus one can estimate $F(g)$ by estimating the normalized characters  values $\frac{\chi(g)}{\chi(1)}$.
A lot has been done on this issue in recent years  and Garion and Shalev used it to prove that almost all elements of $G$ are commutators (a well-known conjecture of Ore asserting that in finite non-abelian simple group every element is a commutator has been proved recently \cite{LOST}).
But one needs more: $g$ should be equal to $[a,b]$ when the pair $\{a,b\}$ generates $G$.
A well-known result of Dixon \cite{Di}, Kantor--Lubotzky \cite{KL} and Liebeck--Shalev \cite{LSh} says:

\begin{thms}\label{thms5.4}
Almost all pairs $(a,b)\in G\times G$ generate $G$.
\end{thms}

Garion and Shalev show that the distribution of the commutators $[a,b]$ over generating pairs $(a,b)$ is approximately the same as over all pairs and  Theorem \ref{theorem5.3} then follows.

\subsection{$G$ simple and $n\geq 3$:}
We saw in Theorem \ref{thms5.4} that almost all pairs of elements of $G$ generate $G$.
This implies that for $n\geq 3$, almost all $n-1$ tuples of generators of $G$ are \emph{redundant}, i.e., a proper subset of the $n$-set already generates $G$.

The following important result was proved by Gilman \cite{Gi} for $n\geq 4$ and extended by Evans \cite{Eva93B} to $n\geq 3$:
\begin{thms}\label{thms5.5}
Let $G$ be a finite simple group and $3\leq n\in \N$.
Then:
\begin{enumerate}[\upshape a.]
\item All the redundant vectors in $\tilde X_n(G)$ lie in the same connected component $Y$ of $\tilde X_n(G)$.
\item The group $\Aut(F_n)$ acts on $\bar Y$, the projection of $Y$ to $\bar E_n(G)$, as the alternating or the symmetric group of degree $|\olin Y|$.
\end{enumerate}
\end{thms}

\begin{cors}\label{corolary5.6}
Let $G$ be a finite simple group.
For $n$ large enough, $\Aut(F_n)$ acts on $\bar E_n(G)$ as the alternating or the symmetric group of degree $|\bar E_n(G)|$.
\end{cors}

As observed by Pak \cite{Pa}, part (a) of the theorem together with Theorem \ref{thms5.4} imply that the graph $\tilde X_n(G)$ has a huge connected component $Y$ whose size is at least $(1-\varepsilon)|\tilde X_n(G)|$ for every $\varepsilon>0$ when $|G|\to \infty$.
So, Wiegold conjecture is essentially true.
For some questions this is enough, but for others (see \S\ref{subsect5.4} below) it is crucial to know that there are no very small connected components in $\tilde X_n(G)$.

Part (b) of the theorem is also very interesting.
We will see in \S\ref{section8} its analogues when $G$ is a compact group.

Part (a) of the theorem is proved using the notion of ``spread''.
(A notion that was first introduced in \cite{BrWi}).

\begin{defns}\label{defns5.6}
A 2-generated group $G$ is said to have spread $r$ if for every non-identity elements $y_1,\ldots,y_r\in G$, there exists $z\in G$ such that $G=\langle y_i,z\rangle$ for every $i=1,\ldots,r$.
\end{defns}

\begin{thms}[Breuer--Guralnick--Kantor \cite{BGK}]\label{thms5.7}
All finite simple groups have spread 2.
\end{thms}

Many (but not all) of them have even spread 3, but we need only 2 in order to prove Theorem \ref{thms5.5}(a) as follows:
Let $z=(z_1,\ldots,z_n)$ and $y=(y_1,\ldots,y_n)$ be two redundant generating vectors.
We can assume $z_i=y_j=1$ for some $i$ and $j$ and after permuting the elements $z_n=y_n=1$.

We can further assume that $z_1\neq 1\neq y_2$.
As $G$ has spread 2, there exists $w\in G$ with $\langle z_1,w\rangle=\langle w,y_2\rangle=G$.
Now, as $\la z_1,\ldots,z_{n-1}\ra=G$ we can move $z=(z_1,\ldots,z_{n-1},1)$ to $(z_1,\ldots, z_{n-1},w)$ and then, using the fact that $\langle z_1,w\rangle=G$, to $(z_1,y_2,\ldots,y_{n-1},w)$.
But we also have $\langle y_2,w\rangle=G$, so the latter can be transformed to $(y_1,y_2,\ldots,y_{n-1},w)$ and finally to $(y_1,\ldots,y_{n-1},1)$ as $\langle y_1,\ldots,y_{n-1}\rangle=G$.

The proof of part (b) of Theorem \ref{thms5.5} is more involved:
One first shows that the action of $\Aut(F_n)$ on $\bar Y$ is double transitive.
Then it is shown that there exists $\beta\in\Aut(F_n)$ acting non-trivially on $\bar Y$ but moves at most $|G|$ elements.
Now, an old result of Bochert (from 1897) asserts that a double transitive permutation subgroup of $\Sym(N)$ with a non-identity element which moves less then $\frac13(N-2\sqrt N)$ elements, must contain $\Alt(N)$.
>From Theorem \ref{thms5.4} we know that $|Y|$ grows like $|G|^{n-1}$ so one can deduce that the action $\Aut(F_n)$ on $\bar Y$ contains $\Alt(\bar Y)$ at least when $G$ is large and $n\geq 4$ (and with little more precision  one sees that this is true for every $G$ and $n\geq 3$).

Another corollary of Theorem \ref{thms5.5}a.\ is that the following is an equivalent reformulation of Wiegold's Conjecture 6.1 (a formulation which is easier to generalize to infinite groups -- see \S\ref{section8} and \S\ref{section9}).

\begin{conj}\label{conjecture5.9}
Let $n\geq 3$ and $\Psi: F_n\twoheadrightarrow G$ an epimorphism onto a finite simple group.
Then $F_n$ has a proper free factor $H$ such that $\Psi(H)=G$.
\end{conj}

Conjecture \ref{conjecture5.9} asserts that for some set of free generators $(g_1,\ldots,g_n)$ of $F_n$, $\varphi(g_1),\ldots,\varphi(g_{n-1})$ generate $G$ and this is the same as saying that the $n$-generating vector $(\varphi(x_1),\ldots,\varphi(x_n))$ is connected to the redundant vector $(\varphi(g_1),\ldots,\varphi(g_n))$.

Let us denote by $E'_n(G)$ the set of all epimorphisms $\Psi:F_n\twoheadrightarrow G$ for which $F_n$ has a proper factor $H$ with $\Psi(H)=G$.
So Wiegold's Conjecture \ref{conjecture5.1} predicts that for finite simple group $G$, $E'_n(G)=E_n(G)$ and $\olin Y=\olin X_n(G)$ for $n\geq 3$.
For the purpose of the product replacement algorithm, the fact that $\olin Y$ is almost all of $\olin X_n(G)$ is just as good.
One starts the algorithm with a redundant vector and so the random walk can take it to almost every other generating vector.
(For the rate of mixing see \cite{Pa}, \cite{LP} and \cite{DS}).
But there are several good reasons to want to know that $\olin Y=\olin X_n(G)$.
This would imply that every presentation of $G$ is equivalent to a minimal one (i.e.\ one with minimal numbers of generators - see \S\ref{section2}) and that the minimal possible of relations can be obtained with a minimal number of generators.
Another application to the representation theory of $\Aut(F_n)$ will be described below, but let us first summarize the very few partial results known toward Wiegold conjecture.

Essentially all the known results are elaborations of the seminal paper of Gilman \cite{Gi}.
Here is the current state of affairs:

\begin{thms}\label{thms5.8}
Let $G$ be a finite simple group and $n\in\N$, then $\tilde X_n(G)$ is connected in the following cases:
\begin{enumerate}[\upshape (i)]
\item \cite{Gi} $G=\PSL_2(p)$, $p$ prime and $n\geq 3$.
\item \cite{Eva93B} $G=S_{z}(2^n)$, the Suzuki groups, or $G=\PSL_2(2^n)$ and $n\geq 3$.
\item \cite{MW} $G=PSL_2(3^p)$, $p$ prime and $n\geq 3$.
\item \cite{Gar08} $G=PSL_2(p^r)$, $p$ prime, $r\in\N$ and $n\geq 4$.
\item \cite{AvGa} $G$ is a finite simple group of Lie rank at most $r$ and $n\geq f(r)$ for a suitable function $f$ depending only on $r$.
\item $G=A_k$, $k\leq 10$ and $n=3$ (see \cite[Theorem 2.5.6.]{Pa}).
\end{enumerate}
\end{thms}

In light of Theorem \ref{thms5.5}a, proving a connectivity result for $\tilde X_n(G)$ amounts to showing that every non-redundant vector $(g)=(g_1,\ldots,g_n)\in {X}_n(G)$ is connected to a redundant one.
We can therefore assume that $g_1,\ldots,g_{n-1}$ generate a proper subgroup of $G$.

The proof of (i) by Gilman is heavily based on the explicit known list of subgroups of $G$ and the same remark is true for the works
of Evans, McCullough--Wanderley and Garion proving (ii), (iii) and (iv) respectively.
For (v), Avni and Garion are using the work of Larsen and Pink \cite{LaPi}
which gives some quantitative description of the possible subgroups of $G$.

Part (v) should be compared with Proposition \ref{proposition3.1} and the result of Nikolov thereafter.
The point here is that the function depends only on $r$ and not on the  defining field.
But there is a price for it, the function $f(r)$ in the proof in \cite{AvGa} grows quite fast (exponentially) with $r$.

Part (vi) was proved using ad-hoc arguments and by help of computer calculations (see \cite{Pa} and the references therein).

So altogether, Wiegold conjecture is known only in very limited cases.
In the next subsection, we will give further motivation to prove the conjecture or even a weak form of it.

\subsection{Representations of $\Aut(F_n)$}\label{subsect5.4}
For many years, it has not been known if $\Aut(F_n)$, the automorphism group of the free group $F_n(n\geq 2)$, or $B_n$ the braid group on $n$ strands $(n\geq 4)$ or $M_g$ - the mapping class group of a closed surface of genus $g\geq 2$ are linear groups, i.e., whether they have faithful linear representations over the field $\C$ of complex numbers.
It has been felt anyway, that all the  problems are similar and a solution to one of them would lead to a solution of all others.
Moreover, it was shown in \cite{DFG} that $B_4$ is linear iff $\Aut(F_2)$ is linear.
In \cite{FP}, Formanek and Processi showed that $\Aut(F_n)$ is \emph{not} linear for $n\geq 3$, leaving the case of $n=2$ open.
The proof was very special for these groups and did not shed any light on $B_n$ (which is a subgroup of $\Aut(F_{n})$).
Moreover,  \cite{BrHe} shows that their method of ``poison subgroup'' can not be applied at all to $M_g$.
On the other hand, Bigelow \cite{Bi} and Krammer \cite{Kr} showed that $B_n$ are all linear and so $\Aut(F_2)$ is linear in spite of the fact that $\Aut(F_n)$, $n\geq 3$ are not.

A new method to produce representations of $\Aut(F_n)$ onto arithmetic groups has been developed recently in \cite{GL09}, but none of these representations is faithful.

We want now to explain a way of looking at this problem which suggests that the difference between $n=2$ and $n\geq 3$ in the linearity question is related to the difference between Wiegold conjecture (Conjecture \ref{conjecture5.1}) for $n\geq 3$ on the connectivity of $\olin X_n(G)$ and $n=2$ where Garion--Shalev showed a strong non-connectivity (Theorem \ref{theorem5.3}).

In fact, even more can be said.
The proof of the non-linearity of $\Aut(F_n)$, $n\geq 3$ in \cite{FP} strongly suggests the following stronger statement:
\begin{conj}\label{conjecture5.11}
Let $n\geq 3$ and $\rho:\Aut(F_n)\to GL_k(\C)$ be a linear representation.
Then $\rho(\Inn(F_n))$ is virtually solvable, where $\Inn(F_n)$ is the group of inner automorphisms of $F_n$.
\end{conj}

Let us now show:
\begin{claim}\label{claim5.10}
If Weigold's conjecture \ref{conjecture5.1} is true, then Conjecture \ref{conjecture5.11} is true when $H=\overline{\rho(\Aut(F_n))}$, the Zariski closure of $\rho(\Aut(F_n))$, is connected.
\end{claim}

As this is only a conditional result, we shall only sketch the proof:

Assume there is such a $\rho$ with $\rho(\Inn F_n)$ not virtually solvable, then by dividing $H$ by its solvable radical we can assume $H$ is semisimple, and even simple, by choosing a suitable factor.
Furthermore, by \cite[Theorem 4.1]{LL} there is a specialization of $\rho$ so that $\rho(\Aut(F_n))$ is in $GL_r(k)$ for some number field $k$ and the Zariski closure $L$ is defined over $k$ but it is still isomorphic to the simple group $H$ over $\C$.

As $\Aut(F_n)$ is finitely generated, $\rho(\Aut(F_n))$ is a subgroup of $L(O_S)$, the $S$-integers of $k$, where $S$ is a finite set of primes of $k$ and $O$ is the ring of integers of $k$.
We can further arrange that $L$ is simply connected and then to apply the strong approximation theorem for linear groups of \cite[Window 9]{LuSe} to deduce that $\rho(\Aut(F_n))$ is almost dense in the congruence completion $L(\hat O_S)$.
The same applies also to $\rho(\Inn(F_n))$, since it is also Zariski dense in $L$.
This implies that for almost every prime ideal $\calP$ in $O_S$, the projection of  $\rho(\Inn(F_n))$ to the finite semisimple group $M=L(O_S/\calP)/Z$ is onto (where $Z$ is the center).
Let $N_\calP$ be the kernel of the map from $\Inn(F_n)$ to $M$.
This subgroup $N_\calP$ is also normal in $\Aut(F_n)$ as it is equal to $\Inn(F_n)$ intersected with the kernel of the map from $\Aut(F_n)$ to $M$.
This means that $N_\calP$ is a characteristic subgroup  of $\Inn(F_n)\simeq F_n$.

Now, $M$ is a product of a bounded number of finite simple groups (a bound independent of $\calP$).
Thus, for infinitely many finite simple groups $G$, when $\Aut(F_n)$  acts on the set of kernels of epimorphisms from $F_n$ onto $G$, it has orbits of a bounded length.
This means that the graph $\olin X_n(G)$ has a component of bounded size in contradiction to its connectivity predicted by Wiegold conjecture.
\qed

The proof shows that much less than Wiegold conjecture is needed.
For example, the following would suffice:
Given a Chevalley group scheme $G$, prove that $X_n(G(\F_q))$ cannot have components of bounded size when $q\to \infty$.
(In fact, using Chebotarev density theorem, it would suffice to assume $q$ is a prime.)

Here is a ``baby version'':
Prove that $F_n$ has no characteristic subgroup $N$ such that $F_n/N$ is a finite simple group.

Another version that would make it:
Prove that for every epimorphism $\Psi: F_n\to G(\F_q)$ (as above), $\Psi(\Phi)$ contains an unbounded number of conjugacy classes  in $G(\F_q)$ when $q\to\infty$.
Here $\Phi$ is the set of primitive elements of $F_n$ (i.e., those belonging to a basis of $F_n$).
Clearly, $\Psi(\Phi)$ is a union of conjugacy classes.
It is not difficult to see that if $\Psi$ corresponds to a redundant generating vector in $\tilde X_n(G(\F_q))$, then $\Psi(\Phi)=G(\F_q)$.
So, again one should only look at non-redundant vectors.

We end this section by remarking that the connection we observed above goes also in the opposite direction:
If $\bar X_n(G(\F_{q_i}))$ has bounded size, say $l$, component for infinitely many primes $q_i$, then we get maps
\[
\rho_{q_i}:\Aut(F_n)\to \Aut(G(\F_{q_i})^l)=\Aut(G(\F_{q_i})^l)\rtimes\Sym (l).
\]
>From these one can cook up a characteristic $0$ representation $\rho$ of $\Aut(F_n)$, with $\rho(\Inn(F_n))$ being Zariski dense in $G(\C)$, i.e., contradicting Conjecture \ref{conjecture5.11}.

\section{Profinite groups}\label{section7}

In this section we discuss the analogous problem in the category of profinite groups.
We will show a strong positive result to all the questions mentioned above, in the context of profinite groups.

The main technical tool which is responsible for it is the following result of Gash\"{u}tz.
Since it is so important, we give an elegant proof due to Roquette (see \cite[Lemma 17.7.2]{FJ}).

\begin{lems}[Gash\"utz Lemma]\label{lemma6.1}
Let $\pi:G\thra H$ be an epimorphism between two finite groups.
Assume $d(G)\leq d$ and let $z_1,\ldots,z_d\in H$ be a set of $d$ elements with $\la z_1,\ldots,z_d\ra=H$.
Then there exist $y_1,\ldots,y_d\in G$ with $\pi(y_i)=z_i$ for $i=1,\ldots,d$ and $G=\la y_1,\ldots,y_d\ra$.
In other wards, any generating $d$-vector of $H$ can be lifted to a generating $d$-vector of $G$.
\end{lems}

\begin{proof}
For every subgroup $B\leq G$ and every $(\ulin{t})=(t_1,\ldots,t_d)\in H$ with $\la t_1,\ldots,t_d\ra=H$, we denote
\[
\varpbt=\#\{(b_1,\ldots,b_d)\in B^d|\pi(b_i)=t_i\text{ for }i=1\ldots,d\text{ and }\la b_1,\ldots,b_d\ra=B\}
\]

\begin{claim*}
The function $\varpbt$ depends only on $B$ and not on $(\ulin{t})$, i.e., it is constant on $(\ulin{t})\in X_d(H)$.
\end{claim*}
This is proved by induction on the size of $B$.
Assume it is true for every proper subgroup $A$ of $B$ and we will prove it for $B$:
Now, if $\pi(B)\lneqq H$ then $\varpbt=0$ for every $(\ulin{t})$ and we are done.
Otherwise,
\[
\varpbt=|K_B|^d-\sum_{A\lneqq B}\varpt{A}
\]
where $K_B=\Ker(\pi:B\thra H)$.
By induction, $\varpt{A}$ is independent of $(\ulin{t})$ and so $\varpbt$ is also independent of $(\ulin{t})$.
This proves the claim.

We deduce now that $\varpt{G}$ is independent of $(\ulin{t})$.
Let $x_1,\ldots,x_d\in G$ be elements such that $\la x_1,\ldots,x_d\ra=G$ (such elements exist since $d(G)\leq d$).
Thus $\varphi_G((\pi(x_i)))>0$, hence also $\varphi_G((z_i))>0$, which is exactly what the lemma says.
\end{proof}

Gash\"utz Lemma does not hold if $G$ is infinite.
When $G=F_d$, its failure is ``measured'' by the number of connected components of $\tilde X_d(H)$.
The lemma has the following corollary for the PRA-graphs:

\begin{cors}\label{newcorollary6.2}
Let $\pi:G\thra H$ be an epimorphism between two finite groups.
Then the induced map $\tilde \pi:\tilde X_n(G)\to \tilde X_n(H)$ is onto for every $n\geq d(G)$.
\end{cors}

Standard inverse limit arguments imply that Lemma \ref{lemma6.1} holds also when $G$ and $H$ are profinite groups and ``generating'' means generating in the topological sense, i.e., generating a dense subgroup.
We can now deduce:
\begin{props}\label{proposition6.2}
Let $F=\hat F_d$ be the free profinite group on $d\in\N$ generators.
If $G$ is a profinite group and $\pi_1,\pi_2$ two epimorphisms from $F$ onto $G$, then there exists a (continuous) automorphism $\alpha$ of $F$, $\alpha\in\Aut(F)$ such that $\pi_1\circ\alpha=\pi_2$.
\end{props}
\begin{proof}
Say $\hat{F}_d=\hat{F}(x_1,\ldots,x_n)$ and denote $z_i=\pi_2(x_i)$ for $i=1,\ldots,d$.
Let $y_1,\ldots,y_d\in \hat{F}_d$ with $\pi_1(y_i)=z_i$ for $i=1,\ldots,d$, and $\overline{\la y_1,\ldots,y_d\ra}=\hat{F}_d$.
Such $y_i$'s exist by Lemma \ref{lemma6.1}.

Let $\alpha$ be the homomorphism from $\hat{F}_d$ to $\hat{F}_d$ sending $x_i$ to $y_i$.
Then, $\alpha$ is an epimorphism and hence an automorphism since every epimorphism from a finitely generated profinite group onto itself is an automorphism.
Moreover, $\pi_1\circ\alpha(x_i)=\pi_1(y_i)=z_i=\pi_2(x_i)$ and so $\pi_1\circ \alpha=\pi_2$ as claimed
\end{proof}

It follows that profinite presentations satisfy all the good properties discussed in \S\ref{section2}.
For example, Waldhausen conjecture:
i.e., if $N\lhd \hat{F}_d$  and $d(\hat{F}_d/N)<d$ then $N$ contains a primitive element of $\hat{F}_d$ (i.e., an element which belongs to a basis of $\hat{F}_d$).

Recall that $\la X; R\ra$ is a profinite presentation for a profinite group $G$, if $R$ is a subset of the free profinite group $\hat{F}_d$ on $X=\{x_1,\ldots,x_d\}$ and $G\simeq \hat{F}_d/\overline{\twola R\twora}$ where $\overline{\twola R \twora}$ is the topological closure of the normal closure of $R$ in $\hat{F}_d$.
We denote by $\hat r(G)$ the minimal possible size of $R$ over all possible profinite presentations of $G$.
Proposition \ref{proposition6.2} now implies that $\hat{r}(G)$ is obtained with a presentation on $d=d(G)$ generators and in all such presentations (as they are all equivalent by the proposition). Moreover, in \cite{L3} it is shown that it is obtained only in these representations.
We mention in passing the long standing open problem:
\begin{prob}\label{problem6.3}
Let $G$ be a finite group.
Is $r(G)=\hat{r}(G)$? where $r(G)$ is the minimal number of relations needed to define $G$ in the discrete category (see \S\ref{section2}) and $\hat r(G)$ the number needed in the profinite category.
\end{prob}
Clearly $\hat{r}(G)\leq r(G)$ but there is no single example of a finite group where a strict inequality is known.
The potential difference between $\hat{r}(G)$ and $r(G)$ was used in \cite{L2} to prove the Mann--Pyber conjecture on the normal subgroup growth of free group.
For $\hat{r}(G)$ one has an exact formula in terms of the cohomology of $G$ (see \cite{GK99} and \cite{L3}) while estimating $r(G)$ is highly non-trivial.
In fact we do not know any lower bound on $r(G)$ for any group which is not, at the same time, also a lower bound for $\hat{r}(G)$.
In \cite{GKKL1} and \cite{GKKL2} presentations of finite simple group are studied and the connections and differences between discrete and profinite presentations are discussed in length.

\section{Compact Lie groups}\label{section8}
In this section $F=F_n$ will denote again the discrete free group on $n$ generators, while $G$ will be a connected compact Lie group.
In this case, for every $n\geq 2$, the set of $n$-tuple $(y_1,\ldots,y_n)\in G^n$, for which $\oline{\la y_1,\ldots,y_n\ra}=G$, is open dense and of full measure in $G^n$.
Thus measurewise the set $G^n$ and the set $\Epi(F,G)=\{\varphi:F\to G|\olin{\varphi(F)}=G\}$ are indistinguishable.
We shall therefore look at the action of $\Aut(F_n)$ on $G^n$.

The main result here is due to Gelander \cite{Ge}
\begin{thms}\label{theorem7.1}
Let $G$ be a compact connected Lie group and let $n\geq 3$.
The action of $\Aut(F_n)$ on $G^n$ in ergodic.
\end{thms}

Theorem \ref{theorem7.1} was conjectured by Goldman \cite{Go2} who proved it for $G=SU(2)$.
He also showed that $n\geq 3$ is necessary (compare Theorem \ref{theorem5.3} above; the reason is similar) and that a proof for the semisimple case would imply the general case.

The space $G^n$ can be thought, also, as the space of $n$-generated marked subgroups of $G$ (i.e., the subgroup is marked by an ordered $n$-generating set of it, so each subgroup appears many times in this space).
A result of the kind of Theorem \ref{theorem7.1} implies that for every measurable property of subgroups of $G$ we have a 0-1 law, i.e., either the property is true for almost all subgroups or it is false for almost all of them, since the subset of marked subgroups with this property is a measurable subset of $G^n$ which is invariant under $\Aut(F_n)$.
(Note that the action of $\Aut (F_n)$ changes the generating set but not the generated subgroup of $G$).
An interesting example of such a property is the \emph{spectral gap property}:
Let $\Gamma\leq G$ be a dense subgroup.
The left translation action of $\Gamma$ on $G$ induces a unitary representation of $\Gamma$ on $\calL^2(G)$.
The complement to the constant functions $\calL_0^2(G)=\{f\in\calL^2(G)|\int f d\mu=0\}$ is $\Gamma$ invariant.
We say that the action of $\Gamma$ on $G$ has a {\itshape spectral gap} if the action of $\Gamma$ on $\calL_0^2(G)$ does not weakly contain the trivial representation.
It is well-known that this happens, for example, if $\Gamma$ has Kazhdan property $(T)$ (see \cite{L1} and the references therein).
It also happens for $G=SU(2)$ with very special choices of $\Gamma$ (see \cite{L1}) based on Deligne's solution to the Ramanujan conjecture.
Altogether, for every semisimple connected compact Lie group $G$, there is such $\Gamma$.
Such $\Gamma$ is responsible for the affirmative answer to Ruziewicz problem (see \cite{L1}).
But it is not known what is the behavior of the generic group with respect to the spectral gap property.
(But see \cite[Theorem 1.4]{LPS}.)
In \cite{Fi}, Fisher pointed out that Theorem \ref{theorem7.1}
implies that either almost all subgroups of $G$ have the spectral
gap or almost all do not. In any event it implies that the set of
$n$-tuples $(n\ge 3)$ in $G^n$ which generate a group with the
spectral gap property is dense in $G^n$. An analogous problem in
the finite groups world is: Essentially all finite simple groups
$G$ have a subset of $k$ generators $\Sigma$ w.r.t.\ which the
Cayley graph $Cay(G;\Sigma)$ is an $\varepsilon$-expanders ($k$
and $\varepsilon$ are independent of $G$), and this is also a
spectral gap property (see \cite{KLN}). But it is not known what
is the behavior of the random set of generators of finite simple
groups, except of the case of the family $\{\PSL_2(p)|p \text{
prime}\}$ where Bourgain and Gamburd \cite{BoGa} showed that
almost all $k$-tuple of elements ($k\geq 2$) give rise to
expanders.

Let us now sketch the proof of Theorem \ref{theorem7.1} for the case $G$ is a semisimple group.
Assume the contrary; let $A\subset G^n$ be an $\Aut(F_n)$ almost invariant measurable subset which is neither null nor conull.
Since $\Aut(F_n)$ is countable we can assume, by replacing $A$ by $\bigcap_{\alpha\in\Aut(F_n)}\alpha(A)$,  that $A$ is $\Aut(F_n)$-invariant.
Now, the action of $G^n$ on itself is clearly ergodic, so at least one of the components, say the first one $G=G_1$, does not preserve $A$.
For $(\ulin{g})=(g_2,\ldots,g_n)\in G^{n-1}$ denote $A_{(\ulin{g})}=\{g\in G| (g,g_2,\ldots,g_n)\in A\}$

\begin{claim*}
For a set of positive measure of $(\underline g)=(g_2,\ldots,g_n)\in G^{n-1}$, the set $A_{(\ulin{g})}$ is neither null nor conull.
\end{claim*}

\begin{proof}
By Fubini, $\mu(A)=\int_{(\ulin{g})} \mu(A_{(\ulin{g})})d\ulin{g}$.
Now, $\mu(A)>0$.
We can throw out those $(\ulin{g})$ with $\mu(A_{(\ulin{g})})=0$ (they contribute measure 0).
So, if for almost all the rest $\mu(A_{(\ulin{g})})=1$, then for every $h\in G_1$, $h\cdot A$ is almost $A$ (since $hA_{(\ulin{g})}\sim A_{(\ulin{g})}$).
But we assumed that $A$ is not $G_1$-almost-invariant.
Thus for a positive measure of $(\ulin{g})$, $0<\mu(A_{(\ulin{g})})<1$ as claimed.

Fix now a point $(\ulin{g})=(g_2,\ldots,g_n)$ in the subset of the claim, such that $\{g_2,g_3\}$ generates a dense subgroup of $G$ (recall that we noticed that the set of such pairs is open, dense and of full measure in $G^2$ -  so such $(\ulin{g})$ does exist!).
The orbits of the action of $\la g_2,g_3\ra$ by left translation on $G_1=G$ coincides with the (projection to the first factor of the) action of the Nielsen moves $\la L(1,2), L(1,3)\ra$ on $\{(g,g_2,g_3,\ldots,g_n)|g\in G\}$ (where $L(1,i)$ sends $(z_1,z_2,\ldots,z_n)$ to $(z_iz_1,z_2,\ldots,z_n)$).
Let $A_1=A_{(\ulin{g})}$ for the $(\ulin{g})$ chosen above.
By our assumption $A_1$ is neither null nor conull.
But on the other hand it is invariant under $\la g_2,g_3\ra$, a dense subgroup of $G$, a contradiction, since every dense subgroup acts ergodically on $G$.
The theorem is now proven.
\end{proof}

We conclude by mentioning another result of Gelander which is proved by similar methods.

\begin{thms}\label{theorem7.2}
Let $n\geq 3$ and $G$ a connected compact Lie group.
Assume $\Gamma\leq G$ is an $(n-1)$-generated dense subgroup.
Then every $n$ elements $s_1,\ldots,s_n\in G$ admit an arbitrary small deformation $t_1,\ldots,t_n$ with $\Gamma=\la t_1,\ldots,t_n\ra$.
In other wards, the set $\{f\in\Hom(F_n,G)|f(F_n)=\Gamma\}$  is dense in $\Hom(F_n,G)$.
\end{thms}

This theorem can be used to prove that given a simple compact Lie
group $G$ containing a dense Kazhdan subgroup, then for some $n$,
any $n$ elements can be $\varepsilon$-deformed (for every
$\varepsilon > 0$) to generate a Kazhdan subgroup of $G$.

\section{Non-compact simple Lie groups}\label{section9}

Let now $G$ be a non-compact simple real Lie group.
In this case one cannot expect to have an ergodic action of $\Aut(F_n)$ on $G^n=\Hom(F_n,G)$.
The representations with discrete image, on one hand and those with dense image on the other hand form two disjoint $\Aut(F_n)$-invariant subsets with non-trivial interior and so the action is not ergodic.
Very little seems to be known about the decomposition of $G^n$ under the $\Aut(F_n)$-action in the general case.
But, recently Minsky has revealed the picture for $G=\PSL_2(\R)$ and $\PSL_2(\C)$. His (somewhat surprising) description shows that this decomposition can be quite delicate, but very interesting.

In the rest of this chapter let $G$ be either $\PSL_2(\R)$ or $\PSL_2(\C)$ and $n\geq 3$.
(For the case $n=2$, the situation is similar to what we show in \S\ref{subsection5.2} for finite groups and \S\ref{section8} for compact groups: the trace of the commutator is an invariant which is preserved by $\Aut(F_n)$ and hence the action is far from being ergodic.)
It will be more convenient to talk about the character variety $X_n(G)=\Hom(F_n,G)/G$.
(We are ignoring the difference between this quotient and the geometric invariant category quotient---see \cite{LM85}---as anyway the representations which are not Zariski dense in $G$ form a measure zero set.)
We will describe how the $\Out(F_n)$ decomposition of $X_n(G)$ following \cite{Mi}.
The reader is referred to that paper and the references thereof, for unexplained  notions and proofs.

Let $\calD$ be the subset of the (equivalent classes) of \emph{faithful discrete representations}.
It contains $\calS$ the Schottky representations.
In fact, it is known (and by no means trivial) that $\calS$ is precisely the interior of $\calD$ and $\olin \calS=\calD$.
The action of $\Out(F_n)$ on $\calS$ is properly discontinuous.

At the other side we have $\calE=\olin E_n(G)$--- the set of representations with dense image, which is an open subset of $X_n(G)$.
The complement of $\calD\cup\calE$ in $X_n(G)$ is the set of all representations which are either discrete but not faithful or are non-discrete but not dense.
This is a measure zero set, so can be ignored for our purpose.
The naive expectation has been that while $\Out(F_n)$ acts properly discontinuously on $\calS$, it would act ergodically on $\calE$ --- a phenomenon that might be seen as an extension of Wiegold's Conjecture \ref{conjecture5.1} for finite simple groups and Gelander's Theorem (Goldman's Conjecture) \ref{theorem7.1} for compact groups.
But this is \emph{not} the case!
In fact Minsky main result in \cite{Mi} is:
\begin{thms}\label{theorem8.1}
There is an open subset of $X_n(G)$, strictly larger than
$\calS$, the set of Schottky representations, which is $\Out(F_n)$
invariant and on which $\Out(F_n)$ acts properly discontunus.
\end{thms}

The set promised in the theorem is $\calP\calS$ the
primitive-stable representations, to be defined below. While it
has a non-empty open intersection with $\calE$, the set of dense
representations, it has an empty intersection with  $R$- the set
of redundant representations, i.e. those representations
$\rho:F_n\to G$, for which there exists a proper free factor $A$
of $F_n$ with $\rho(A)$ dense in $G$. (Compare to Theorem
\ref{thms5.5} and Conjecture \ref{conjecture5.9}. Note that for
$G$ compact $R$ is conull in $G^n$ - see \S\ref{section8}.) The
set $R$ is open. One is tempted to suggest:\pagebreak
\begin{conj}\label{conjecture8.2}\
\begin{enumerate}[\upshape (a)]
\item The action of $\Out(F_n)$ on $R$ is ergodic.
\item $R\cup \calP\calS$ is conull in $X(G)$.
\end{enumerate}
\end{conj}
If true, this conjecture gives a nice satisfactory picture:
$X_n(G)$ is, up to a set of measure zero, a union of two
$\Out(F_n)$-invariant open subsets $\calP\calS$ and $R$. On the
first $\Out(F_n)$ acts properly discontinuous and on the second it
acts ergodically. But, at this point this is just wishful
thinking. (See a remark added in proof at the end of this
section).

Let us now define $\calP\calS$ and describe Minsky's main ingredients.

Let $C$ denote the Cayley graph of $F_n$ with respect to the free generators $x_1,\ldots,x_n$, and $\partial C=\partial F_n$ the boundary of $F_n$, i.e.\ the rays from an initial vertex to infinity on the graph.
Let  $\partial ^2 F_n=(\partial F_n\times \partial F_n)\setminus \Delta$ where $\Delta$ is the diagonal.
Thus $\partial^2 F_n$ is the set of biinfinite (oriented) lines on $C$.

To each $w\neq1$ in $F_n$ we associate a biinfinite line, i.e.\ a point $\olin w=(\infty_1,\infty_2)$ in $\partial ^2F_n$, i.e.\ a biinfinite word obtained by concatenating infinitely many copies of a representative of $w$.
If $g\in F_n$, then the point of $\partial ^2F_n$ associated with $gwg^{-1}$ is $(g\infty_1,g\infty_2)$.
We denote by $\olin{\olin w}$ the $F_n$-orbit of $\olin w$.
Let $\calP$ be the  subset of $\partial ^2F_n$ of all points associated with primitive elements of $F_n$.
It is clearly invariant under the action of $F_n$, and $\Out(F_n)$ acts on the set $\calB$ of $F_n$-orbits.

A representation $\rho:F_n\to G$ and a base point $x_0$ in the symmetric space $\bbH$ (which is either $\bbH ^2$ if $G=\PSL_2(\R)$ or $\bbH^3$ if $G=\PSL_2(\C)$) gives rise to a unique map $\tau_{\rho,x_0}\to\bbH$ mapping the origin of $C$ to $x_0$, which is $\rho$-equivariant and maps each edge to a geodesic).
Every element of $\calB$ is represented by an $F_n$-invariant set of infinite lines which is  mapped to a family of broken geodesic paths in $\bbH$.

\begin{defns}
A representation $\rho:F_n\to G$ is \textit{primitive-stable} if there are constants $K,\delta$ in $\R_+$ and a basepoint $x_0\in\bbH$ such that $\tau_{\rho,x}$ takes the lines representing the primitive elements $\calP$ to $(K,\delta)$-quasi geodesic.
This means that for some $K,\delta\in \R_+$, for any two vertices $v_1, v_2$ on a line in $\calP$, $\frac1K \dist_{\bbH}(v_1,v_2)-\delta\leq \dist_C(v_1,v_2)\leq K\dist_{\bbH}(v_1,v_2)+\delta$.
\end{defns}

If there is one such basepoint, then any basepoint will do, at the expense of increasing $\delta$.
The set of primitive stable representatives is $\Aut(F_n)$-invariant.
Its image in $X_n(G)$ will be denoted $\calP\calS$.

Schottky representations give rise to quasi-isometric embeddings of $C$ in $\bbH$, so $\calS\subset\calP\calS$.
The converse is not true, but Minsky showed that if $\rho\in\calP\calS$, then for every proper free factor $A$ of $F_n$, $\rho(A)$ is Schottky group.
The set $\calP\calS$ like $\calS$ is also open; and the action of $\Out(F_n)$ is properly discontinuous.

The crucial point in proving this last claim is that the image of the set
$\{\alpha\in\Aut(F_n)|\|\alpha(w)\|\leq c\|w\|\;\forall\text{ primitive }w\}$ in $\Out(F_n)$ is finite, where $c$ is any finite constant and for $g\in F_n$ we denote by $\|g\|$ the length of its cyclically reduces word (i.e., the minimal length in its conjugacy class).

The above facts are relatively simple to deduce from the basic definitions.
The nontrivial fact is that $\calP\calS$ is indeed larger than $\calS$.
To this end Minsky shows that one representation $\rho_0$ at the boundary of $\calS$ (i.e., $\rho_0$ in $\calD\setminus\calS$) is primitive-stable (in fact, he gives a method to produce many such examples, but one suffices!).
Since $\calP\calS$ is open it implies that some open neighborhood of $\rho_0$ is also in $\calP\calS$ - such a neighborhood has a nontrivial open intersection with $\calE$. So indirectly one deduces the existence of many primitive-stable representations with dense images--- even though only a discrete one is explicitly constructed!

To construct the discrete non-Schottky primitive-stable representation, Minsky appeals to a result of Whitehead which, using what nowadays is called the Whitehead graph, gives a necessary criterion for a word in $F_n$ to be primitive.
Given $g\in F_n$ define the graph $Wh(g)$ to be the graph with $2n$ vertices denoted by the generators $\{x_i\}_{i=1}^n$ and their inverses $\{x_i^{-1}\}_{i=1}^n$.
A pair $(a,b)$ of vertices is an edge if $ab^{-1}$ appears in $g$ or in a cyclic permutation of $g$ (which is the same as saying it appears in $g$ or $g$ starts with $b^{-1}$ and  end with $a$) call this last edge the additional edge if it does not appear anyway in $g$.

\begin{thms}[Whitehead]\label{thms9.4}
Let $g$ be a cyclically reduced primitive element in $F_n$.
Then by eliminating (at most) one vertex, $Wh(g)$ becomes a non-connected graph.
\end{thms}

Whitehead's result gives a simple sufficient criterion for a word $w$ to be ``blocking'' --- i.e., one which cannot appear as a subword of any cyclically reduced primitive element.
This is the case if $Wh(w)$, minus the additional edge, contains a cycle which passes through all the vertices of the graph.
It is easy to see that this is the case for $\beta^2=([x_1,x_2][x_3,x_4]\cdots[x_{2m-1},x_{2m}])^2$ as an element of $F_n$, $n=2m$.
It also follows now that $\beta^2$ is not inside any proper free factor of $F_n$.

Let now $\Sigma$ be a surface of genus $m$ with one boundary which is a curve represented by $\beta$ in $\pi_1(\Sigma)=F_{2m}$.
Let $\rho:\pi_1(\Sigma)\to\PSL_2(\R)$ be a discrete representation with $\rho(\beta)$ being parabolic.
It is well-known that such $\rho$ exists in this case.
It is a very special case of a general result asserting  that every simple curve $\gamma$ on the boundary of $3$-dimensional handlebody gives rise to a geometrically finite representation into $\PSL_2(\C)$ for which $\rho(\gamma)$ is parabolic (see \cite{Mi} and the references therein).
As $\rho(F_{2m})$ contains a parabolic element, it is not Schottky.

Minsky then proves that this $\rho$ is primitive stable. This is
done as follows: Let $Y$ be the convex hull of the limit set of
$\rho(\pi_1(\Sigma))$. In our case, as $\rho(\pi_1(\Sigma))$ is a
non-uniform lattice in $\PSL_2(\R)$, $Y$ is actually equal to
$\bbH^2$, but this is not crucial for the general case. Let
$Z=Y/\rho(\pi_1(\Sigma))$--- the convex core of $\rho$. This is a
surface with a unique cusp. Minsky shows that all the primitive
elements of $F_n=\pi_1(\Sigma)$ are represented by geodesics in a
fixed compact set $K\subset Z$. The idea is that in order to leave
a compact set, a primitive element must wind around the cusp and
this is prohibited by the blocking property deduced from
Whitehead's Lemma. The existence of this compact $K$ implies the
quasi-isometric condition for primitive elements, in a way similar
to the standard argument that a group acting cocompactly is
quasi-isometric to the space upon which it acts.

This finishes the sketch of the proof of Theorem \ref{theorem8.1} for $n$ even and some modifications give the general case.

The theorem leaves various interesting problems. Define
$\calP\calS'(G)$ to be the set of all (equivalent classes of)
representations of $F_n$ where restrictions to proper free factors
are Schottky. So $\calP\calS\subset \calP\calS'$ and Minsky shows
that this is a proper inclusion. He asks whether $\calP\calS$ is
the interior of $\calP\calS'$. He also shows that no point outside
$\calP\calS'$ can be in the domain of discontinuity of $\Out(F_n)$
acting on $X_n(G)$. Thus a positive answer to this question will
show that $\calP\calS$ is exactly the domain of discontinuity  for
the action of $\Out(F_n)$ on $X_n(G)$. Together with Conjecture
\ref{conjecture8.2} this will give a nice picture of the action of
$\Out(F_n)$ on $X_n(G)$ for these two cases of $G$. One can
speculate to suggest that similar picture holds also for
$G=\PSO(r,1)$ for $r\geq 2$. For these $G$'s, at least the
definitions make sense. We do not know what even to expect for the
situation to be for higher rank simple Lie groups $G$. The work of
Minsky shows that the naive extension of Weigold--Goldman
Conjectures (8.1 and 6.1) is false. But it still seems somewhat
likely that the action of $\Aut(F_n)$ on $R(G)$---the set of
redundant representations is always ergodic. This will be a
beautiful analogue of the theorem of Gilman and Evans (Theorem
\ref{thms5.5}) and Gelander's Theorem (Theorem \ref{theorem7.1}).

We end this section by describing a recent work of Glasner \cite{Gla} that shows that this is indeed the case for two families of simple locally compact groups.
So let us now switch notations and assume that $G$ is either $\PSL_2(K)$ where $K$ is a non-archimedean local field of characteristic $0$ or $G=\Aut^+(T_k)$---the group of orientation preserving automorphisms of the $k$-regular tree $T_k$, $k\geq3$.
(It is a simple group of index two in the full automorphism group of $T_k$.)
Note also that $\PSL_2(K)$ is acting on a tree; the Bruhat--Tits tree associated with it.

The Schottky subgroups of these $G$'s were studied in detail in \cite{Lub91}.
It is shown there that the subset $\calS$ of the Schottky representations is an open and closed subset of $H_n(G)=G^n$ and of $X_n(G)=\Hom(F_n,G)/G$.
The action of $\Out(F_n)$ on $\Hom(F_n,G)/G$ is not studied there, but from the discussion it is not difficult to see that $\Out(F_n)$ acts properly discontinuously on $\calS$.
Let now $\olin E_n(G)$ denote the subset of $X_n(G)$ of all the dense representations.
\begin{thms}[Glasner \cite{Gla}]\label{thms9.5}
Let $G$ be either $\PSL_2(K)$ or $\Aut^+(T_k)$.
Then for every $n\geq 3$, $\Out(F_n)$ acts ergodically on $\olin E_n(G)$.
\end{thms}

In fact, he shows that $\Aut(F_n)$ acts ergodically on the set of all dense representations in $\Hom(F_n,G)=G^n$.
Before sketching the proof, let us first mention that for these $G$'s, $\calS\cup\olin E_n(G)$ is far from covering the whole space.
We also have an open subset of all the representations of $F_n$ whose image lie in the compact open subgroup (the stabilizer of a vertex).

Glasner's proof is based on two main ingredients. The first is a
result of Weidman \cite{Wei02} asserting that if $\rho(F_n)$ is
dense in $G$ then $\rho(w)$ is elliptic  (i.e.\ fixes a  vertex)
for some \emph{primitive} element $w$ of $F_n$. This implies that
every $n$-tuple in $E_n(G)$ is conjugate mod $\Aut(F_n)$ to an
$n$-tuple of type $(w,g_2,\ldots,g_n)$ with $\rho(w)$ elliptic.
Glasner shows further that $\rho(g_2)$ can be made to be
hyperbolic. Then he uses another result (proved in \cite{AbGl} for
$\Aut^+(T_k)$ and in \cite{Gla} for $\PSL_2(K)$): for almost every
elliptic element $a$ and almost every hyperbolic element $b$ the
group generated by $a$ and $b$ is dense in $G$. From this, he
applies some arguments of a similar nature to Gelander proof of
Theorem \ref{theorem7.1},  to deduce the theorem. Along the way he
shows that $R(G)$ is conull in $E_n(G)$.

All these results of Minsky and Glasner seem to indicate that only
the tip of the iceberg has been revealed. It looks like a rich and
interesting theory should be explored here for general non-compact
Lie groups (or other locally compact groups).

\subsection*{Added in proof} Conjecture 9.2(a) has been proved
recently by Gelander and Minsky \cite{GeMi}.  In fact they proved
it for every simple $k$-group defined over a characteristic $0$
local field $k$.  Their work explains the difference between
$G=PSL_2(\R)$ or $PSL_2(\C)$ for which the action of $\Out(F_n)$
on $\bar E_n(G)$ is not ergodic (see Theorem 9.1) and the group
$G=PSL_2 (\Q_p)$ for which it is ergodic (Theorem 9.5). The
crucial difference is that for the latter, almost every dense
representation of $F_n(n\ge 3)$ to $G=PSL_2(\Q_p)$ is redundant (a
fact whose proof by Glasner uses Weidman \cite{Wei02} in a crucial
way).

\section{The mapping class group action on surface group representations}\label{section10}

In the previous sections we studied the action of $\Aut(F_n)$ on $\Hom(F_n,G)$ ( and of $\Out(F_n)$ on $\Hom(F_n,G)/G$) for  various groups $G$.
In principle, one can do this not only for $F_n$ but also for any finitely generated group $\Gamma$.
A case of special interest is $\Gamma=\Pi_g$ ---  the fundamental group of a closed surface $\Sigma_g$ of genus $g\geq 2$.
Indeed, this case has been studied in the literature in great detail as it is related to classical geometric and topological topics such as Fricke--Teichm\"uller spaces.
A comprehensive survey is given by Goldman \cite{Go1}, who is responsible, to a large extent, for the modern systematic development of the theory.
In this section we mention only few points out of this theory.
Our main goal is to call the attention to a particular direction which is not covered in \cite{Go1}; the study of the action of $\Out(\Pi_g)$, the mapping class group, on $\Epi(\Pi_g,G)/G$ when $G$ is a finite group.
This issue came out in a recent paper of Dunfield and Thurston \cite{DT} where finite sheeted covers  of random $3$-manifolds are studied.
It suggests developing a theory of the kind described in Sections  \ref{section5} and \ref{section6}, for $\Pi_g$ instead of $F_n$.
One may, for example, suggest an analogous conjecture to Weigold's, a proof of which (or even of a weak form of it) would imply that the mapping class groups are not linear.

But let us start with $G$ being infinite:
It is of interest to note that for $\Gamma=F_n$ the study of the $\Aut(\Gamma)$ action on $\Hom(\Gamma,G)$ has started with $G$ finite in presentation theory, as described in the Sections \ref{section2}--\ref{section6}, and only later a systematic study for $G$ compact or semisimple has emerged.
On the other hand for $\Gamma=\Pi_g$ the most classical case is the study of $\Hom(\Pi_g,\PSL_2(\R))/\PSL_2(\R)$.
The faithful discrete representations form a connected component which is exactly the space classifying the equivalent classes of conformal structures on $\Sigma_g$, or also equivalence classes of hyperbolic structures on $\Sigma_g$.
But there are more components which are indexed by the Euler class $e:\Hom(\Pi_g,\PSL_2(\R))/\PSL_2(\R)\to H^2(\Sigma_g;\Z)\simeq\Z$ whose image is $\{2-2g,\ldots,2g-2\}$, i.e., $4g-3$ connected components.
The components $e^{-1}(\pm(2-2g))$ are two copies of the Teichm\"uller space which differ by the choice  of orientation.
On these two, $M_g=\Out(\Pi_g)$, which is classically known as the mapping class group of $\Sigma_g$, is acting properly discontinuously and a lot of study has been devoted to this action by many authors (see \cite{Go1} and the references therein).
Much less is known about the action on the other $4g-5$ components.
Goldman conjectures that $M_g$ acts ergodically on each of these.
If $\PSL_2(\R)$ is replaced by a connected compact Lie group,then it was indeed proved by Pickrell and Xia \cite{PX} that $M_g$ acts ergodically on every component of $\Hom(\Pi,G)/G$.
The special case $G=SU(2)$ was proved by Goldman \cite{Go3} who conjectured the general case for $\Gamma=\Pi_g$ as well as for $\Gamma=F_n$ as discussed in Section \ref{section8}.
If $G$ is semisimple compact group, then the number of the connected components of $\Hom(\Pi_g,G)/G$ is equal to the order of the fundamental group of $G$.
The same applies for complex semisimple groups $G$, but is not true in general, for example
$\Hom(\Pi_g,\widetilde{SL_3(\R)})/\widetilde{SL_3(\R)}$ is not connected.

A wealth of additional information is given in \cite{Go1}, but we will move now to the case when $G$ is a finite group, which is not discussed there.

In \cite{DT}, Dunfield and Thurston suggest an interesting model to produce random $3$-manifolds.
It briefly goes like that:
It is well-known that every closed $3$-manifold $M$ has an Heegard splitting, i.e., it can be presented as a union of two handlebodies of genus $g$, $H_1$ and $H_2$ which are glued along their boundaries, each of which is a genus $g$ closed surface.
Their idea is to use this as a way to produce closed $3$-manifolds of Heegard genus (at most) $g$ in the following way:
Fix $g$ and fix a set of generators $S$ for the mapping class group $M_g$ of $\Sigma_g$.
Take a random walk along the Cayley graph of $M_g$ with respect to $S$.
This will produce a random element $\varphi\in M_g$.
Use this random $\varphi$ to glue the boundary of $H_1$ ---  an handlebody of genus $g$ --- to a copy of it, $H_2$, along the boundary.
This will give the resulting ``random'' $3$-manifold.
They were interested in finite covers of such random manifolds and in questions of the following type:
Given a finite group $G$, what is the probability that a random $3$-manifold $M$ of genus $g$ as above, has a finite sheeted cover $M'$ with a cover group isomorphic to $G$?
This is really the question:
What is the probability that there is an epimorphism from $\pi_1(M)$ onto $G$?
Now, $\pi_1(M)$ can be described in the following way:
Start with
$\pi_1(\Sigma_g)=\la a_1,b_1,\ldots,a_g,b_g|\prod_{i=1}^{g}[a_i,b_i]=1\ra$
the fundamental group of the surface $\Sigma_g$.
Gluing $H_1$ to it ``kills'' $a_1,\ldots,a_g$ and we get the free group on $b_1,\ldots,b_g$.
Then gluing $H_2$, the second copy, amounts to dividing $\pi_1(\Sigma_g)$ further by $\varphi(a_1),\ldots,\varphi(a_g)$ to get $\pi_1(M)$.
(Note that $\varphi\in M_g$ gives an element of $\Aut\pi_1(\Sigma_g)$ which is well-defined only up to inner automorphism but the normal closure of $\varphi(a_1),\ldots,\varphi(a_g)$ is well-defined.)

Now, let $\rho:\Pi_g\thra G$ be an epimorphism, then it ``survives'' in the above process if and only if $a_i$ and $\varphi(a_i)$ are in $\Ker\rho$ for every $i=1,\ldots,g$.
There are many such epimorphisms $\rho$ (a very good estimate is given in \cite{LiSh2}, at least for the most interesting case, when $G$ is a finite simple group).
We can ask the above question in a different way now:
Start with an epimorphism $\rho:\Pi_g\to G$ with $\rho(a_i)=1$ for every $i=1,\ldots,g$.
What is the probability that for a random $\varphi\in M_g$, $\varphi^{-1}(\Ker\rho)$ still contain $a_1,\ldots,a_g$?
For the discussion of this question and the interesting answer(s) we refer the reader to \cite{DT}.
For our context, what is relevant is the steps taken in \cite{DT} to study the action of $M_g$ on the set of all kernels of epimorphisms from $\Pi_g$ onto $G$.

This last action cannot be expected to be transitive in general.
In fact, the epimorphism $\rho:\Pi_g\thra G$ induces a map  $H_2(\Pi_g,\Z)\to H_2(G;\Z)$ and thus to every kernel $\Sigma_g\thra G$ one associates an invariant $[c]\in H_2(G,\Z)/\Out(G)$.
By using a ``stabilization'' result of Livingston \cite{Li}, it is shown in \cite{DT}:

\begin{thm*}
Let $G$ be a non-abelian finite simple group.
Then for all sufficiently large $g$, the orbits of $\Epi(\Pi_g,G)$ under $M_g=\Out(\Pi_g)$ correspond bijectively to $H_2(G,\Z)/\Out(G)$.
Moreover, the action of $M_g$ on each orbit is by the full alternating group of that orbit.
\end{thm*}

This theorem is the analogue of Corollary \ref{corolary5.6}.
It will be of interest to give a quantitative estimate of the $g$ needed for a given $G$, as in Corollary \ref{newcorollary3.3}.
It will be even more remarkable if one can prove a ``Wiegold's Conjecture'' in this context, i.e.\ that for $g\geq 3$ (and actually maybe even $g\geq 2$) $M_g$ acts transitively on all the kernels of $\Pi_g\thra G$ with the same invariant in $H_2(G,\Z)/\Out(G)$.
One can then imitate the discussion in \S\ref{subsect5.4} (and just like there, a weaker statement suffices: there are no bounded size orbits) to deduce that $\Aut(\Pi_g)$ is not linear.
>From this last statement one can conclude that $M_{g+1}$ is not linear.
As of now, this is a long standing open problem.

\newcommand{\etalchar}[1]{$^{#1}$}
\providecommand{\bysame}{\leavevmode\hbox to3em{\hrulefill}\thinspace}
\renewcommand{\MR}{\relax\ifhmode\unskip\space\fi }
\providecommand{\MRhref}[2]{%
  \href{http://www.ams.org/mathscinet-getitem?mr=#1}{#2}
}
\providecommand{\href}[2]{#2}

\end{document}